\documentclass{article}
\usepackage{amsfonts}
\usepackage{amsmath}
\usepackage{amssymb}
\usepackage{latexsym}
\newtheorem{thm}{\sc Theorem}[section]

\newtheorem{lem}[thm]{\bf Lemma}
\newtheorem{prop}[thm]{\bf Proposition}
\newtheorem{defn}[thm]{\bf Definition}
\newtheorem{rem}[thm]{\bf Remark}
\newtheorem{exmp}[thm]{\bf Example}

\newcommand{\field}[1]{\mathbb{#1}}

\newcommand{\Q }{\field{Q}}
\newcommand{\Z }{\field{Z}}

\def\proof{{\parindent0pt {\it Proof.\ }}}

\def\w.dim{{\rm w.dim}}

\def\max{{\rm max}}
\def\qf{{\rm qf}}
\def\qt{{\rm Q}}

\def\pd{{\rm pd}}
\def\fd{{\rm fd}}

\begin{document}
\baselineskip=15pt

\begin{center}
{\small Comm. Algebra 32 (10) (2004) 3937-3953}
\end{center}
\vspace{1cm}

\begin{center}
{\Large\bf Trivial Extensions Defined by
\\Coherent-like Conditions}

\bigskip\bigskip
{\large\bf Salah-Eddine Kabbaj $^{1,}$\footnote{The first author
was supported by KFUPM under project \# MS/TRIVIAL/224.\\
Correspondence: Salah-Eddine Kabbaj, Department of Mathematics,
KFUPM, P.O. Box 5046, Dhahran 31261, Saudi Arabia; E-mail:
kabbaj@kfupm.edu.sa} and Najib Mahdou $^2$}

\bigskip
$^1$ Department of Mathematics, King Fahd University of Petroleum
and Minerals\\Dhahran, Saudi Arabia

$^2$ Department of Mathematics, FST Fez-Sa\" iss, University S. M.
Ben Abdellah\\Fez, Morocco
\end{center}

\bigskip\bigskip
\noindent{\large\bf Abstract.} This paper investigates
coherent-like conditions and related properties that a trivial
extension $R:= A\propto E$ might inherit from the ring $A$ for
some classes of modules $E$. It captures previous results dealing
primarily with coherence, and also establishes satisfactory
analogues of well-known coherence-like results on pullback
constructions. Our results generate new families of examples of
rings (with zerodivisors) subject to a given coherent-like
condition.

\bigskip
\begin{section}{INTRODUCTION}

All rings considered in this paper are commutative with identity
elements and all modules are unital. Let $A$ be a ring and $E$ an
$A$-module. The trivial ring extension of $A$ by $E$ is the ring
$R:= A \propto~E$ whose underlying group is $A \times E$ with
multiplication given by $(a, e)(a', e') = (aa', ae'+a'e)$.
Considerable work, part of it summarized in Glaz's book \cite{Gz2}
and Huckaba's book (where $R$ is called the idealization of $E$ in
$A$) \cite{H}, has been concerned with trivial ring extensions.
These have proven to be useful in solving many open problems and
conjectures for various contexts in (commutative and
non-commutative) ring theory. See for instance
\cite{DS,Fo,KM,Le,M1,M2,PR,Po,R}.

A ring $R$ is coherent if every finitely generated ideal of $R$ is
finitely presented; equivalently, if $(0:a)$ and $I\cap J$ are
finitely generated for every $a\in R$ and any two finitely
generated ideals $I$ and $J$ of $R$ \cite{Gz2}. Examples of
coherent rings are Noetherian rings, Boolean algebras, von Neumann
regular rings, valuation rings, and Pr\"ufer/semihereditary rings.
The concept of coherence first sprang up  from the study of
coherent sheaves in algebraic geometry, and then developed, under
the influence of Noetherian ring theory and homology, towards a
full-fledged topic in algebra. During the past 30 years, several
(commutative) coherent-like notions grew out of coherence  such as
finite conductor, quasi-coherent, $v$-coherent, $n$-coherent, and
-to some extent- GCD and G-GCD rings (see the respective
definitions in the beginning of Sections 2 and 3). Noteworthy is
that both the ring-theoretic and homological aspects of coherence
run through most of these generalizations (see for instance
\cite{Gz1}).

This paper investigates coherent-like conditions and related
properties that a trivial extension $R:= A\propto E$ might inherit
from the ring $A$ for some classes of modules $E$. It captures
previous results dealing primarily with coherence \cite{Gz2,R},
and also establishes satisfactory analogues of well-known
coherence-like results on pullback constructions \cite{GH} (see
also \cite{BR,DP,DKM,DKMS,FG}). Our results generate new families
of examples of rings (with zerodivisors) subject to a given
coherent-like condition.

The second section provides a ring-theoretic
approach. We first extend the definition of a $v$-coherent domain to
rings with zerodivisors and develop a theory of these rings parallel
to Glaz's study of finite conductor, quasi-coherent, and G-GCD rings
\cite{Gz1}. Afterwards, we study the possible transfer of all these
notions for various trivial extension contexts. Thereby, new examples are
provided which, particularly, enrich the current literature with new
classes of  coherent-like rings with zerodivisors.

The third section treats the homological aspect. We first study
conditions under which trivial extensions yield (strong)
$n$-coherent rings \cite{C,CK,DKM,DKMS}. Due to reciprocal effects
\cite[Section 2]{C}, we also deal with the $(n,d)$-rings of Costa,
i.e., those in which $n$-presented modules \cite{Bo2} have
projective dimension at most $d$. In particular, the second part
of this section is devoted to Costa's second conjecture that one
may characterize the $(n,d)$-property intrinsically by
ideal-theoretic-conditions \cite{C}. We explore the scope of the
validity of this conjecture in various trivial extension
non-coherent contexts. Recall at this point that Costa's second
conjecture is valid in the class of coherent rings \cite{CK}. This
fact was behind our motivation for studying large classes of
coherent-like rings. The paper closes with an independent result
showing that this conjecture holds in the class of finite
conductor domains (resp., rings) for $n\leq2$ and $d=1$ (resp.,
$d=0$). The general case is still elusively open.

The following diagram of commutative rings summarizes the
relations between the coherent-like notions involved in this
paper:\bigskip

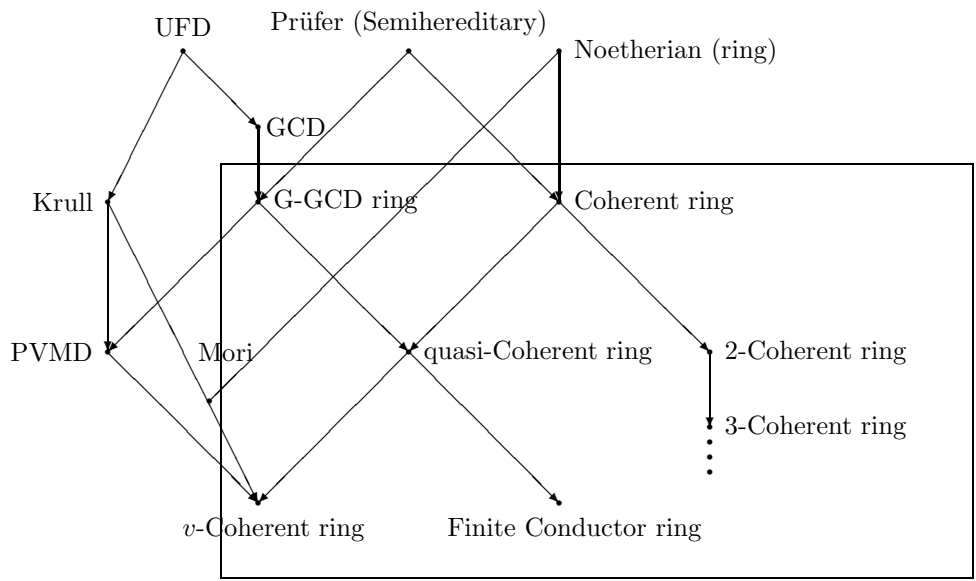
\[\setlength{\unitlength}{1mm}
\begin{picture}(100,65)(10,-20)
\put(10,20){\line(1,-2){10}}
\put(10,20){\vector(0,-2){20}}
\put(20,40){\vector(-1,-2){10}}
\put(20,40){\vector(1,-1){10}}
\put(30,30){\vector(0,-1){10}}
\put(30,20){\vector(-1,-1){20}}
\put(10,0){\vector(1,-1){20}}
\put(20,0){\vector(1,-2){10}}
\put(70,40){\vector(0,-1){20}}
\put(70,20){\vector(-1,-1){20}}
\put(70,20){\vector(1,-1){20}}
\put(30,20){\vector(1,-1){20}}
\put(50,0){\vector(1,-1){20}}
\put(50,0){\vector(-1,-1){20}}
\put(90,0){\vector(0,-1){10}}
\put(50,40){\vector(-1,-1){20}}
\put(50,40){\vector(1,-1){20}}
\put(70,40){\line(-1,-1){46.5}}

\put(25,25){\line(10,0){100}}
\put(25,-30){\line(10,0){100}}
\put(25,25){\line(0,-5){55}}
\put(125,25){\line(0,-5){55}}

\put(20,40){\circle*{.7}}
\put(20,42){\makebox(0,0)[b]{UFD}}
\put(10,20){\circle*{.7}}
\put(8,20){\makebox(0,0)[r]{Krull}}
\put(10,0){\circle*{.7}}
\put(8,0){\makebox(0,0)[r]{PVMD}}
\put(23.5,-6.5){\circle*{.7}}
\put(22,0){\makebox(0,0)[l]{Mori}}
\put(30,20){\circle*{.7}}
\put(32,20){\makebox(0,0)[l]{G-GCD ring}}
\put(30,30){\circle*{.7}}
\put(31,30){\makebox(0,0)[l]{GCD}}
\put(70,40){\circle*{.7}}
\put(72,40){\makebox(0,0)[l]{Noetherian (ring)}}
\put(70,20){\circle*{.7}}
\put(72,20){\makebox(0,0)[l]{Coherent ring}}
\put(50,0){\circle*{.7}}
\put(52,0){\makebox(0,0)[l]{quasi-Coherent ring}}
\put(90,0){\circle*{.7}}
\put(92,0){\makebox(0,0)[l]{$2$-Coherent ring}}
\put(70,-20){\circle*{.7}}
\put(72,-22){\makebox(0,0)[t]{Finite Conductor ring}}
\put(30,-20){\circle*{.7}}
\put(32,-22){\makebox(0,0)[t]{$v$-Coherent ring}}
\put(90,-10){\circle*{.7}}
\put(92,-10){\makebox(0,0)[l]{$3$-Coherent ring}}
\put(90,-12){\circle*{.7}}
\put(90,-14){\circle*{.7}}
\put(90,-16){\circle*{.7}}
\put(50,40){\circle*{.7}}
\put(50,42){\makebox(0,0)[b]{Pr\"ufer (Semihereditary)}}
\end{picture}\]
\end{section}

\bigskip

\begin{section}{RING-THEORETIC APPROACH}

A ring $R$ is  quasi-coherent (resp., finite conductor) if
$(0:a)$ and $a{_1}R\cap ... \cap a{_n}R$ (resp., $bR\cap cR$) are
finitely generated ideals of $R$ for any finite set of elements $a$
and $a{_1}, ..., a{_n}$ (resp., $b,c$) of $R$
\cite{BAD,Gz1,Z}. Also, $R$ is called a G-GCD ring if every
principal ideal of $R$ is projective and the intersection of any two
finitely generated flat ideals of $R$ is a finitely generated flat
ideal of $R$ \cite{AA,Gz1}.

\begin{subsection}{$v$-Coherent Rings with Zerodivisors}

In view of Glaz's recent work on finite conductor, quasi-coherent, and
G-GCD rings \cite{Gz1}, we first extend the definition of a
$v$-coherent domain \cite{FG,GH,N1,N2} to rings with  zerodivisors. For this
purpose, we review some terminology related to basic operations on
fractional ideals in an arbitrary ring (i.e., not necessarily a domain). Let
$R$ be a commutative ring and let $\qt(R)$ denote the total ring of
quotients of $R$.  By an ideal of $R$ we mean an integral ideal of $R$. Let
$I$ and $J$ be two  nonzero
fractional ideals of $R$. We define the fractional ideal $(I:J)=\{x\in
\qt(R)\mid xJ\subset  I\}$. We denote $(R:I)$ by $I^{-1}$ and
$(I^{-1})^{-1}$ by $I_v$ (called the $v$-closure of $I$). A nonzero
fractional ideal $I$ is said to be invertible if $II^{-1}=R$, divisorial (or
a $v$-ideal) if $I_v=I$, and $v$-finite if $I_v=J_v$ (or, equivalently, if
$I^{-1}=J^{-1}$) for some finitely generated fractional ideal  $J$ of $A$.
The $v$-operation on $R$ is not necessarily a $*$-operation, since, in
general,  $(a)_{v}\neq (a)$ when $a$ is a zerodivisor of $R$. However, the
other basic properties of the  $v$-operation on integral domains
\cite[(32.1)(2)\&(3), (32.2)(a)\&(b)]{Gi} carry over to arbitrary rings.
\begin{defn} \label{2.1}\rm
A ring $R$ is $v$-coherent if $(0:a)$ and $\bigcap_{1\leq i\leq n}Ra_{i}$
are  $v$-finite ideals of $R$ for any finite set of elements $a$ and $a{_1},
..., a{_n}$ of $R$.
\end{defn}

\begin{prop} \label{2.2}
Let $R$ be a ring and let's consider the following assertions:\\
{\rm(1)} $I^{-1}$ is $v$-finite for any finitely generated ideal $I$ of
$R$.\\
{\rm(2)} $I{_v}\cap J_v$ is $v$-finite for any two finitely generated ideals
$I$ and $J$ of $R$.\\
{\rm(3)} $\bigcap_{1\leq i\leq n}Ra_{i}$ is $v$-finite for any finite set of
elements $a{_1}, ..., a{_n}$ of $R$. \\
Then {\rm(1)}$\Longrightarrow${\rm(2)}. Moreover, if $R$ is an integral
domain, then the three assertions are equivalent, each of which
characterizes $v$-coherence. $\Box$
\end{prop}
\proof Assume that (1) is true and let $I$ and $J$ be any finitely
generated ideals of $R$. Then there exist two finitely generated
ideals $I_{1}$ and $J_{1}$ such that $I_v = I_{1}^{-1}$ and $J_v =
J_{1}^{-1}$. So,  $I^{-1} = (I_{1})_{v}$ and $J^{-1} =
(J_{1})_{v}$, hence $I_v \cap J_v =(I^{-1} + J^{-1})^{-1}
=((I_{1})_{v} + (J_{1})_{v})^{-1} =((I_{1})_{v})^{-1} \cap
((J_{1})_{v})^{-1} =(I_{1})^{-1} \cap (J_{1})^{-1} =(I_{1} +
J_{1})^{-1}$ which is $v$-finite by hypothesis since $I_{1} +
J_{1}$ is a finitely generated ideal of $R$.

Now, assume that $R$ is an integral domain. Then
(1)$\Longleftrightarrow$(2) is handled by \cite[Proposition 3.6]{FG}, and
(1)$\Longleftrightarrow$(3) always holds since $(\sum_{i=1}^{n}Ra_{i})^{-1}
=\bigcap_{1\leq i\leq n}Ra_{i}^{-1}$ for each $a_i \in R$ and any
integer $n\geq 1$. $\Box$
\bigskip

Clearly, quasi-coherent rings are $v$-coherent, and if $R$ is a
domain, the above definition matches the definition of a $v$-coherent
domain. It is worth recalling that $v$-coherent domains offer a large
context of validity for the so-called Nagata's theorem for the class group
\cite{Ga}. Also, recall from \cite{N1} that PVMDs \cite{Gi,Z} and Mori
domains \cite{BG} are $v$-coherent. Moreover, non-Krull integrally closed
Mori domains \cite{B}  are ($v$-coherent but) not finite conductor \cite{Z}.

Let $(R_{j})_{1\leq j\leq m}$ be a family of rings and $R
=\prod_{j=1}^{m}R_{j}$. For $C=(c_{j})$ and $A{_1}=(a_{1j}), ...,
A{_n}$ $=(a_{nj})\in R$, we have $(0:C) =
\prod_{j=1}^{m}(0:c_{j})$ and $\bigcap_{1\leq i\leq
n}RA_{i}=\prod_{j=1}^{m}\big(\bigcap_{1\leq
  i\leq n}Ra_{i}\big)$. Further, for any ideals $I
=\prod_{j=1}^{m}I_{j}$ and $J =\prod_{j=1}^{m}J_{j}$  of $R$, we have
$I \cap J =\prod_{j=1}^{m}(I_{j} \cap J_{j})$ and $I^{-1}
=\prod_{j=1}^{m}I_{j}^{-1}$. Then $\prod_{j=1}^{m}R_{j}$ is
$v$-coherent if and only if so is $R_j$ for each $j=1, \ldots
,m$. Thus, finite products (for instance, of Mori domains) may provide
us with original examples of $v$-coherent rings with zerodivisors.\bigskip

Let's now examine $v$-coherence for rings of small weak dimension.
Recall first that rings of weak dimension $0$ are precisely the
von Neumann regular rings. Moreover, Glaz showed that for a ring
$R$ of weak dimension $1$ the finite conductor property,
quasi-coherence, and coherence deflate to the mere fact that
$(0:c)$ is finitely generated for every $c\in R$ \cite[Proposition
  2.2]{Gz1}. She also proved that the finite conductor and
quasi-coherence properties coincide for rings of weak dimension
$2$ \cite[Theorem 2.3]{Gz1}. In contrast with these results, the
next example denies any similar effect to weak dimension on
$v$-coherence.

\begin{exmp} \label{2.3}\rm
Let $E$ be a countable direct sum of copies of $\Z/2\Z$ with addition and
multiplication defined component wise, where $\Z$ is the ring of integers.
Let
$R=\Z \times E$ with multiplication defined by
$(a,e)(b,f) =(ab,af+be+ef)$. Then:\\
(1) $\w.dim(R) =1.$\\
(2) $R$ is not coherent.\\
(3) $R$ is a $v$-coherent ring.
\end{exmp}
\proof (1) That $\w.dim(R) =1$ this is handled in \cite[Example 1.3, page
10]{V}.\\
(2) Let $x=(2,0) \in R$. Then $(0:x) =\{(a,e) \in R | (a,e)(2,0) =0\} =
\{(a,e) \in R | (2a,0) =0\} =0 \times E$ which is not a finitely
generated ideal of $R$. Therefore, $R$ is not a coherent ring.\\
(3) Notice first that an element $s\in R$ is regular if and only
if $s=(a,0)$ with $a\in\Z\setminus2\Z$. This easily follows from
the four basic facts: $E$ is Boolean; $2E=0$; $ae=e$ for any
$a\in\Z\setminus2\Z$ and $e\in E$; and for any $e\not=0 \in E$,
there exists $f\not=0 \in E$ such that $ef =0$.

Next, we wish to show that each ideal of $R$ is $v$-finite which
implies that $R$ is $v$-coherent. Let $J$ be an ideal of $R$ and
let $I =\{a \in \Z /(a,e) \in J$ for some $e \in E\}$. Assume $I
=0$. Let $s$ be any regular element  of $R$. Clearly,
$(0,e)=s(0,e)$ for any $e\in E$. It follows that $sJ=J$ and hence
$J^{-1} =Q(R) =(R(0,e))^{-1}$ for any $e\not=0 \in E$. Now, assume
$I=x\Z$, where $x$ is a nonzero integer. We claim that $J^{-1}
=(R(x,0))^{-1}$. Indeed, let $y/s \in Q(R)$, where $y =(a,e) \in
R$ and $s=(b,0)$ is a regular element. It can easily be seen that
$sR =b\Z \times E$. Then $y/s \in J^{-1} \Leftrightarrow yJ
\subseteq sR \Leftrightarrow (a,e)J \subseteq b\Z \times E
\Leftrightarrow aI \subseteq b\Z \Leftrightarrow ax \in b\Z
\Leftrightarrow (a,e)(R(x,0)) \subseteq sR \Leftrightarrow y/s \in
(R(x,0))^{-1}$. Thus, in both cases, $J$ is $v$-finite, as
asserted. $\Box$
\bigskip

While a ring $R$ which is a total ring of quotients trivially is
$v$-coherent, $R$ need not be finite conductor \cite[Example
3.5]{Gz1}. The following construction may offer new contexts that
illustrate this fact.

\begin{exmp} \label{2.4}\rm
Let $(R,M)$ be any local ring with $M^2 =0$. Then:\\
(1)  $R$ is a $v$-coherent ring that is not G-GCD. \\
(2) The following conditions are equivalent:\\
(i) $R$ is a coherent ring;\\
(ii) $R$ is a quasi-coherent ring;\\
(iii) $R$ is a finite conductor ring;\\
(iv) $(0:c)$ is finitely generated for every $c \in R$;\\
(v) $M$ is finitely generated.
\end{exmp}
\proof (1) That $R$ is $v$-coherent this is straightforward since
$R=Q(R)$ is a total ring of quotients. Let $c\not=0\in M$. Then Ann$(c)
=(0:c)=M$. Hence $Rc$ is not
projective (since not free), so that $R$ is not a G-GCD ring
\cite{Gz1}. \\
(2) Clearly, we only need prove $(v)\Longrightarrow(i)$. Assume that $M$ is
finitely generated and let $I$ be a finitely generated proper ideal of $R$.
Let $\{x_{1},\ldots ,x_{n}\}$ be a minimal generating set of $I$ and
consider the exact sequence of $R$-modules:
$$0 \rightarrow  Ker(u) \rightarrow  R^n  \buildrel u \over\rightarrow
I \rightarrow 0$$
where $u(a_{1},\ldots ,a_{n}) = \sum _{i=1}^{n}a_{i}x_{i}$. We claim
that $Ker(u) =\prod M=:M^n$.  Indeed, $M^n \subseteq Ker(u)$ is clear
since $M^2 =0$ and $x_{i} \in M$ for each $i =1, \ldots ,n$. On the
other hand, $Ker(u) \subseteq M^n$ since $\{x_{1},\ldots
,x_{n}\}$ is minimal. Therefore, $Ker(u) =M^n$
is a finitely generated $R$-module (since $M$ is). Hence, $I$ is
finitely presented and thus $R$ is coherent. $\Box$
\end{subsection}

\begin{subsection}{Results of Transfer and Examples}

This subsection  investigates the possible transfer of the coherence
properties for various trivial extension contexts. Our results generate new
families of examples subject to a given coherent-like condition.

For the convenience of the reader, we next  discuss  some basic
facts connected to trivial ring extensions. These will be used
frequently in the sequel without explicit mention.  Let $A$ be a
ring and $E$ an $A$-module and let $R:=A\propto~E$ be the trivial
ring extension of $A$ by $E$. An ideal of $R$ of the form
$I\propto IE$, where $I$ is an ideal of $A$, is finitely generated
if and only if $I$ is finitely generated \cite[page 141]{Gz2}.
Also recall that $R$ has always its Krull dimension equal to the
Krull dimension of $A$ \cite[Theorem 25.1(3)]{H}.

For a general description of modules over a trivial ring
extension, we refer the reader to \cite[pages 140 \& 141]{Gz2}.
Here, we describe a specific type of $R$-modules that play a
crucial role within the $R$-module structure, namely, finitely
generated free $R$-modules and their $R$-submodules. Let $n$ be a
positive integer. Define the ``multiplication'' on $E$ by elements
of $A^n$ within $E^n$ through the natural $A$-bilinear map
$\varphi:A^{n}\times E\longrightarrow E^{n}$  defined by
$ae=\varphi(a,e):=(a_{i}e)_{1\leq i\leq n}$, for any
$a=(a_{i})_{1\leq i\leq n}\in A^{n}$ and $e\in E$. Now let $U$ be
an $A$-submodule of $A^n$ and $E'$ an $A$-submodule of $E^n$ such
that $UE\subseteq E'$. Let $U\propto E'$ denote the set $U\times
E'$ with natural addition and scalar multiplication defined by
\((a, e)(u,e')= (au,ae'+ue)\). Clearly, $U\propto E'$ is an
$R$-module; and, under this notation, the  finitely generated free
$R$-module $R^n$ identifies with $A^n\propto E^n$. Also, $U\propto
E'$ is a finitely generated $R$-module only if $U$ is a finitely
generated  $A$-module. Conversely, let  $M$ be an $R$-submodule of
$R^n$. Set $U:=\{u\in A^{n} | (u,e')\in M\ \textup{for some}\
e'\in E^{n}\}$ and $E':=\{e'\in E^{n} | (u,e')\in M\ \textup{for
some}\ u\in A^{n}\}$. It is easily seen that $U$ and $E'$ are
$A$-modules such that $M \subseteq U\propto E'$. The next example
illustrates the fact that equality does not hold in general.

\begin{exmp} \label{2.5}\rm
Let $(A,M)$ be a local domain which is not a field, $E :=A/M$, and
$R :=A \propto E$ be the trivial ring extension of $A$ by $E$. Let
$J =R(x,1)$, where $x\not=0 \in M$. Set $I =\{a \in A | (a,e) \in
J$ for some $e \in E\}$ and $E' =\{e \in E | (a,e) \in J$ for some
$a \in A\}$. Then $J \subsetneqq I \propto E'$.
\end{exmp}
\proof One may easily check that $I =Ax$ and $E'=E$. Further, we
claim that $(x,0)\in I\propto E'\setminus J$. Deny. We have $(x,0)
=(a,e)(x,1)$ for some $(a,e) \in R$ so that $x =ax$. Hence $a=1
\in M$, the desired contradiction. $\Box$
\bigskip

Nevertheless, it is easily seen that $M =U \propto E'$ if and only
if $0 \propto E' \subseteq M$ if and only if $U \propto 0
\subseteq M$.\bigskip

Example~\ref{2.5} shows that \cite[Theorem 25.1(1)]{H} is not
true. This was confirmed by the author of \cite{H} through a
private e-communication.

\begin{thm} \label{2.6}
Let $(A,M)$ be a local ring and $E$ an $A$-module with $ME =0$. Let $R :=A
\propto E$ be the
trivial ring extension of $A$ by $E$. Then:    \\
{\rm(1)}  $R$ is a $v$-coherent ring that is not G-GCD.\\
{\rm(2)} $R$ is coherent (resp., quasi-coherent, finite conductor)
if and only if $A$ is coherent (resp., quasi-coherent, finite
conductor), $M$ is finitely generated, and $E$ is an
$(A/M)$-vector space of finite rank.
\end{thm}

Before proving Theorem~\ref{2.6}, we establish the following Lemma.

\begin{lem}  \label{2.7}
Under the hypotheses of Theorem~\ref{2.6}, $(0:c)$
is a finitely generated ideal of $R$ for each $c \in R$ if and only if
$(0:a)$
is a finitely generated ideal of $A$ for each $a \in A$,  $M$ is finitely
generated, and $E$ is an $(A/M)$-vector space of finite rank.
\end{lem}
\proof Assume that $(0:c)$ is a finitely generated ideal of $R$ for
each $c \in R$. Let $a \in A$ and set $c :=(a,0) \in R$. Then $(0:c)
=(0:a) \propto E'$, where $E' =\{e \in E | ae=0\}$. Therefore, $(0:a)$
is a finitely generated ideal of $A$. Let
$e\not=0\in E$ and set $c :=(0,e) \in R$. Similar arguments show that
$M$ is a finitely generated ideal of $A$ since $(0:c) =M \propto
E$. Further, assume that $M \propto E= \sum _{i=1}^{n}R(x_{i},e_{i})$, where
$x_{i}\in M$ and $e_{i}\in E$ for each $i=1,..., n$. Then $E \subseteq \sum
_{i=1}^{n}(A/M) e_{i}$, and hence $E$ is an $(A/M)$-vector space of finite
rank.

Conversely, let $c :=(a,e)\not=0 \in R$. If $a$ is invertible in $A$,
then $c$ is invertible in $R$. Then, without loss of generality, we
may assume that $a \in M$. Hence, $(0:c) =\{(b,f) \in R |(ab,be) =0\}
=\{(b,f) \in M \propto E |ab=0\}$ (since if $b$ is invertible in $A$,
then $(b,f)$ is invertible in $R$, and so $c =0$, a contradiction). It
can easily be seen that if
$a =0$ then $(0:c) =M \propto E$, and if  $a \not= 0$ then $(0:c)
=(0:a) \propto E$. In both cases, $(0:c)$ is a finitely generated
ideal of $R$ since $M$ and $(0:a)$ are finitely generated
ideals of $A$ and $E$ is an $(A/M)$-vector space of finite rank, completing
the proof of Lemma~\ref{2.7}.  $\Box$ \bigskip

{\noindent \it Proof of Theorem~\ref{2.6}.} (1) One may easily
verify that $R$ is local with maximal ideal $M \propto E$ and that
each element of $R$ is either a unit or a zerodivisor. Then
$R=Q(R)$ is $v$-coherent. Let $c\not=0\in M$ and $e\not=0\in E$.
Clearly, $(c,e)$ is a zerodivisor. Hence $R(c,e)$ is not
projective (since not free), so that $R$ is not a G-GCD ring
\cite{Gz1}.\\
(2) Assume that $R$ is a coherent ring. By Lemma~\ref{2.7}, it
remains to show that $A$ is coherent.  Let  $I =\sum _{i=1}^n
Aa_{i}$, where $a_{i} \in M$ and set $J
:=\sum_{i=1}^{n}R(a_{i},0)$. Consider the exact sequence of
$R$-modules:
$$0 \rightarrow  Ker(u) \rightarrow  R^n= A^{n}\propto E^n  \buildrel u
\over\rightarrow
J \rightarrow 0$$
where $u((b_{i},e_{i})_{1\leq i\leq n}) =
\sum_{i=1}^{n}(b_{i},e_{i})(a_{i},0)=(\sum_{i=1}^{n}a_{i}b_{i},0)$ since
$a_{i} \in M$ for each $i=1, \ldots ,n$. On the other hand, consider the
exact sequence of $A$-modules:
$$0 \rightarrow  Ker(v) \rightarrow  A^n  \buildrel v \over\rightarrow  I
\rightarrow 0$$ where $v((b_{i})_{1\leq i\leq n}) =
\sum_{i=1}^{n}a_{i}b_{i}$.  Then, $Ker(u) =Ker(v) \propto E^n$.
But $J$ is finitely presented since $R$ is coherent, so $Ker(u)$
is a finitely generated $R$-module and hence $Ker(v)$ is a
finitely generated $A$-module. Therefore, $I$ is a finitely
presented ideal of $A$, so $A$ is coherent.

Conversely, let $J$ be a finitely generated ideal of $R$ and let $S
:=\{(a_{i},e_{i})\}_{1\leq i\leq n}$ be a minimal generating set of $J$,
where $a_{i} \in M$ and $e_{i} \in E$. Consider the exact sequence of
$R$-modules:
$$0 \rightarrow  Ker(u) \rightarrow  R^n  \buildrel u \over\rightarrow  J
\rightarrow 0$$
where $u((b_{i},f_{i})_{1\leq i\leq n})
=\sum_{i=1}^{n}(a_{i},e_{i})(b_{i},f_{i}) =
(\sum_{i=1}^{n}a_{i}b_{i},\sum_{i=1}^{n}b_{i}e_{i})$ since $a_{i} \in
M$ for each $i=1, \ldots ,n$. Further, the minimality of $S$ yields $Ker(u)
=
\{(b_{i},f_{i})_{1\leq i\leq n} \in R^{n} | \sum_{i=1}^{n}a_{i}b_{i}
=0\}$. Let $I :=\sum_{i=1}^{n}Aa_{i}$ and consider the surjective $A$-module
homomorphism $v$ defined above. Then $Ker(v)$ is a finitely generated
$A$-module since $A$ is coherent. Consequently, $Ker(u) =Ker(v) \propto E^n$
is a finitely generated $R$-module. Therefore, $J$ is finitely presented and
hence $R$ is coherent.

Now, assume that $R$ is quasi-coherent. We only need show that
$\bigcap_{1\leq i\leq n}Ra_{i}$ is a finitely generated ideal of $A$
for each $a_{i} \in M$. This is straightforward since
$\bigcap_{1\leq i\leq n}R(a_{i},0) =(\bigcap_{1\leq i\leq n}Aa_{i})
\propto 0$ is a finitely generated ideal of $R$.

Conversely, let $J =\bigcap_{1\leq i\leq n}R(a_{i},e_{i})$, where
$a_{i} \in M$ and $e_{i} \in E$. We may suppose that $J\subsetneqq
R(a_{i},e_{i})$ for each $i=1, \ldots ,n$. Let $(a,e) \in J$. Then,
there exist $b_{i} \in A$ and $f_{i} \in E$ such that $(a,e)
=(b_{i},f_{i})(a_{i},e_{i}) =(a_{i}b_{i},b_{i}e_{i})$ for each $i=1,
\ldots ,n$. We claim that $b_{i} \in M$ for each $i=1, \ldots
,n$. Deny. Then, there exists $j$ such that $b_j$ is invertible in
$A$ and so is $(b_{j},f_{j})$ in $R$. Hence $(a_{j},e_{j})
=(b_{j},f_{j})^{-1}(a,e) \in J$ yielding $J =R(a_{j},e_{j})$, a
contradiction.  Therefore, $(a,e) =(a_{i}b_{i},0)$. It follows that
$J =\bigcap_{1\leq i\leq n}(Aa_{i} \propto 0)=(\bigcap_{1\leq i\leq
n}Aa_{i})\propto 0$ is finitely generated in $R$ since $\bigcap_{1\leq i\leq
n}Aa_{i}$ is by hypothesis finitely generated in $A$. Thus $R$ is
quasi-coherent.

Finally, similar arguments as above with $n =2$ lead to the conclusion
for the finite conductor property, to complete the proof of
Theorem~\ref{2.6}. $\Box$
\bigskip

Next, we explore a different context, namely, the trivial ring
extension of a domain by its quotient field.

\begin{thm}  \label{2.8}
Let $A$ be a domain which is not a field, $K =\qf(A)$, and $R:=A \propto K$
be the trivial ring extension of
$A$ by $K$. Then:\\
{\rm(1)} $R$ is not a finite conductor ring. In particular, $R$ is neither
quasi-coherent nor coherent.\\
{\rm(2)} $R$ is a $v$-coherent ring if and only if $A$ is
$v$-coherent.
\end{thm}
\proof (1) Let $x :=(0,1) \in R$. Then $(0:x)=0\propto K$ is not a
finitely generated ideal of $R$. Therefore, $R$ is not a finite
conductor ring, as asserted.\\
(2) Observe first that $(a,e) \in R$ is regular if and only if $a
\not= 0$, and that $(0:c)=0\propto K$ for any $c:=(0,e)\not=0\in
R$. Further, \cite[Theorem 25.1(4)]{H} yields $\bigcap_{1\leq
i\leq
  n}R(a_{i},e_{i})=\bigcap_{1\leq i\leq n}(Ra_{i}\propto K)=(\bigcap_{1\leq
i\leq n}Ra_{i})\propto K$, for every finite set of elements
$(a_{i},e_{i})_{1\leq i\leq n}$ of $R$ (with $a_{i}\not=0$ for
each $i$). Now, Let $I$ be any nonzero ideal of $A$ and $E$ any
$A$-submodule of $K$ with $IK\subseteq E$ and let $J:=I \propto
E$. By \cite[Theorem 25.10]{H}, we have \(J^{-1} =(I \propto
E)^{-1}= I^{-1}\propto K =(I \propto IK)^{-1}\).  Finally, since
$I\propto IK$ is finitely generated if $I$ is, the conclusion is
straightforward. $\Box$ \bigskip

New examples of original coherent-like rings with zerodivisors with
arbitrary Krull dimensions may stem from Theorems \ref{2.6} \&
\ref{2.8}, as shown by the following constructions.

\begin{exmp}\label{2.9}\rm
Let $K$ be any field and $X_{1},X_{2},...$ be indeterminates over $K$.
Let $n$ be any integer $\geq1$, $A =K[[X_{1},...,X_{n}]]$ the
power series ring in $n$ variables over $K$, and $R :=A \propto
K$. Then, by Theorem~\ref{2.6}, $R$ is  an $n$-dimensional coherent ring
that is not G-GCD. $\Box$
\end{exmp}

\begin{exmp}\label{2.10}\rm
Let $A$ be as in the above example and $R :=A \propto K[Y]$, where
  $Y$ is another indeterminate over $K$. Then, by Theorem~\ref{2.6}, $R$ is
an $n$-dimensional   $v$-coherent ring that is not finite conductor. $\Box$
\end{exmp}

\begin{exmp}\label{2.11}\rm
Let  $R:=\Z[X_{1},...,X_{n-1}] \propto \Q(X_{1},...,X_{n-1})$, where $n$
is any integer $\geq1$, $\Z$ the ring of integers, and $\Q$ the field
of rational numbers. Then, by Theorem~\ref{2.8}, $R$ is an $n$-dimensional
$v$-coherent ring that is not finite conductor. $\Box$
\end{exmp}

More examples are provided in the next section, denying any
possible interplay between some of these coherent-like conditions
and $n$-coherence.

\end{subsection}
\end{section}
\begin{section}{HOMOLOGICAL APPROACH}

For a nonnegative integer $n$, an $R$-module $E$ is $n$-presented
if there is an exact sequence $F_n\rightarrow F_{n-1} \rightarrow
...\rightarrow F_0 \rightarrow E\rightarrow 0$ in which each $F_i$
is a finitely generated free $R$-module  \cite{Bo2}. In
particular, ``$0$-presented'' means finitely generated and
``$1$-presented'' means finitely presented. Throughout, $\pd_R(E)$
and $\fd_R(E)$ will denote the projective dimension and the flat
dimension of $E$ as an $R$-module, respectively. Also $\w.dim(R)$
will denote the weak dimension of $R$.

In 1994, Costa introduced a doubly filtered set of classes of rings, called
the $(n,d)$-rings,  with the aim of obtaining a good understanding of the
structures of some non-Noetherian rings \cite{C}. The Noetherianity forces
the regularity of the $(n, d)$-rings. However, outside Noetherian settings,
the richness of this classification resides in its ability to unify classic
concepts such as von Neumann regular, hereditary, Dedekind, semihereditary,
and Pr\"ufer  rings.

Given nonnegative integers $n$ and $d$, a ring $R$ is called an
$(n,d)$-ring if every $n$-presented $R$-module has projective
dimension $\leq d$ ; and a weak $(n,d)$-ring  if every
$n$-presented cyclic $R$-module has projective dimension $\leq d$
(equivalently, if every $(n-1)$-presented ideal of $R$ has
projective dimension $\leq d-1$). For instance, the
$(0,1)$-domains are the Dedekind domains, the $(1,1)$-domains are
the Pr\"ufer domains, and the $(1,0)$-rings are the von Neumann
regular rings \cite{C}. Costa's paper concludes with a number of
open problems, including his second conjecture that the $(n,d)$-
and weak $(n,d)$-properties are equivalent. This conjecture is
valid in the class of coherent rings \cite{CK}.

The first part of this section studies the transfer of the
$n-$coherence properties (see definitions below) to trivial ring
extensions. Due to reciprocal effects \cite[Section 2]{C}, results
of transfer for the $(n,d)$-property are also provided.  Our
purpose, in the second part, is mainly to test the validity of
Costa's second conjecture for non-coherent contexts.

\begin{subsection}{$n$-Coherence and Strong $n$-Coherence}

Recall from \cite{DKM,DKMS}, for $n\geq 1$, that $R$ is
$n$-coherent if each $(n-1)$-presented ideal of $R$ is
$n$-presented; and that $R$ is strong $n$-coherent if each
$n$-presented $R$-module is $(n+1)$-presented (This terminology is
not the same as that of \cite{C}, where Costa's ``$n$-coherence"
is our ``strong $n$-coherence"). In particular, ``$1$-coherence''
coincides with coherence, and one may view ``$0$-coherence'' as
Noetherianity. Any strong $n$-coherent ring is $n$-coherent, and
the converse holds for $n=1$ or for coherent rings
\cite[Proposition 3.3]{DKMS}. Strong $n$-coherence arose naturally
in Costa's study \cite{C} of the $(n,d)$-rings. As a matter of
fact, every $(n,d)$-ring is strong $\max\{n,d\}$-coherent
\cite[Theorem 2.2]{C}; and an $(n,d)$-ring is strong $r$-coherent
($r<n$) only if it is an $(r,d)$-ring \cite[Theorem 2.4]{C}.
\bigskip

Our main result examines the context of trivial ring extensions of
domains by their respective quotient fields.

\begin{thm} \label{3.1}
Let $A$ be a domain which is not a field, $K =\qf(A)$, and $R:=A
\propto K$ be the trivial ring extension of $A$ by $K$. Let
$n\geq2$ and $d\geq1$ be integers. Then the
following hold:  \\
{\rm(1)} $R$ is not coherent.  \\
{\rm(2)} $R$ is strong $n$-coherent (resp., $n$-coherent) if and only if so is $A$. \\
{\rm(3)} $R$ is an $(n,d)$-ring (resp., a weak $(n,d)$-ring) if
and only if so is $A$.
\end{thm}

The proof of this theorem lies mainly on the following two lemmas
which characterize, respectively, finitely generated and
$n$-presented $R$-submodules of free $R$-modules.\bigskip

Let us fix the notation for the next two results. Let $R$ be as in
Theorem 3.1 and let $H$ be an $R$-submodule of $R^m$, where $m$ is
an arbitrary positive integer. Set $U =\{x \in A^m / (x,e) \in H$
for some $e \in K^{m}\}$ and $E =\{e \in K^m / (x,e) \in H$ for
some $x \in A^{m}\}$.

\begin{lem}  \label{3.2}
Under the above notation, the following statements are equivalent:\\
{\rm(i)} $H$ is finitely generated and $E$ is a $K$-vector space;\\
{\rm(ii)} $U$ is finitely generated and $H =U \propto KU$.
\end{lem}
\proof $i)\Longrightarrow ii)$ Let $H
=\sum_{i=1}^{p}R(x_{i},e_{i})$ ($\subseteq U \propto E$), where
$p$ is a positive integer, $x_{i} \in A^{m}$, and $e_{i} \in
K^{m}$ for each $i =1, \ldots ,p$. Necessarily, $U
=\sum_{i=1}^{p}Ax_{i}$ and $E =\sum_{i=1}^{p}Ae_{i} + KU$. Next
assume $KU \subsetneqq E$. Then there exists a nonzero $K$-vector
space $F$ with finite rank such that $F \oplus KU =E$. Write
$e_{i} =y_{i} + z_{i}$, where $y_{i} \in F$ and $z_{i} \in KU$ for
each $i =1, \ldots ,p$. From above, it follows that $E
=\sum_{i=1}^{p}Ay_{i} \oplus KU$ and thus $F
=\sum_{i=1}^{p}Ay_{i}$. Consequently, $F$ (and hence $K$) is a
finitely generated $A$-module, the desired contradiction. Hence,
$E =KU$. Now let $y \in E =KU$. Then $y
=\sum_{i=1}^{p}b_{i}x_{i}$, where $b_{i} \in K$ for each $i =1,
\ldots ,p$. So $(0,y) =\sum_{i=1}^{p}(0,b_{i})(x_{i},e_{i}) \in
H$. It follows that $0 \propto E \subseteq H$; equivalently, $H =U
\propto KU$.\\
$ii)\Longrightarrow i)$ Straightforward. $\Box$

\begin{lem}\label{3.3}
Let $n$ be an integer $\geq1$. Under the above notation, the following statements are equivalent:\\
{\rm(i)} $H$ is $n$-presented;\\
{\rm(ii)} $U$ is $n$-presented and $H =U \propto KU$.
\end{lem}
\proof $i)\Longrightarrow ii)$ By induction on $n$. Assume $n=1$.
As above, write $H =\sum_{i=1}^{p}R(x_{i},f_{i})$, where $p$ is a
positive integer, $x_{i} \in A^m$, and $f_{i} \in K^m$ for each
$i=1,...,p$. We have $U =\sum _{i=1}^{p}Ax_{i}$ and $E = (\sum
_{i=1}^{p}Af_{i}) + KU$. Let $F'$ be a $K$-vector space such that
$F'\oplus KU =K^m$. Then, $f_{i} =g_{i} + h_{i}$, where $g_{i} \in
F'$ and $h_{i} \in KU$ for each $i=1, \ldots ,p$. Hence $E = (\sum
_{i=1}^{p}Ag_{i}) \oplus KU$. Further, it easily can be seen that
$H =\sum_{i=1}^{p}R(x_{i},g_{i})$. Consider the exact sequence of
$R$-modules:
$$0 \rightarrow Ker(w) \rightarrow R^p \buildrel w \over\rightarrow  H
\rightarrow 0$$ where $w((a_{i},e_{i})_{1, \ldots, p})
=\sum_{i=1}^{p}(a_{i},e_{i})(x_{i},g_{i}) =( \sum
_{i=1}^{p}a_{i}x_{i}, \sum _{i=1}^{p}a_{i}g_{i} + \sum
_{i=1}^{p}e_{i}x_{i})$. Set $W :=\{(a_{i})_{i=1, \ldots ,p} \in
A^p / \sum_{i=1}^{p}a_{i}x_{i} =0$ and $\sum_{i=1}^{p}a_{i}g_{i}
=0\}$ and $E' :=\{(e_{i})_{i=1, \ldots ,p} \in K^p /
\sum_{i=1}^{p}e_{i}x_{i} =0\}$. Clearly, $E'$ is a $K$-vector
space and $Ker(w)$ is a finitely generated $R$-submodule of $R^p$
(since $H$ is finitely presented). By Lemma 3.2, $W$ is finitely
generated and $Ker(w) =W \propto KW$. Moreover, let $(a_{i})_{i=1,
\ldots ,p}\not=0 \in A^p$ such that $\sum_{i=1}^{p}a_{i}x_{i} =0.$
Then, $(a_{i})_{i=1, \ldots ,p} \in E' =KW.$ So there exists
$z\not=0 \in K$ and $(b_{i})_{i=1, \ldots ,p} \in W$ such that
$(a_{i})_{i=1, \ldots ,p} =z(b_{i})_{i=1, \ldots ,p}$. Hence
$\sum_{i=1}^{p}a_{i}g_{i} =z\sum_{i=1}^{p}b_{i}g_{i}=0$, whence
$(a_{i})_{i=1, \ldots ,p} \in W$. Consequently, $W=\{(a_{i})_{i=1,
\ldots ,p} \in A^p / \sum_{i=1}^{p}a_{i}x_{i} =0\}$. Therefore,
the exact sequence of $A$-modules of natural homomorphisms
$$ 0 \rightarrow W \rightarrow A^p \rightarrow  U \rightarrow 0,$$
upon tensoring by the flat $A$-module $R$, yields the exact
sequence of $R$-modules
$$ 0 \rightarrow W \otimes_{A}R \cong W \propto KW =Ker(w) \rightarrow R^p
\rightarrow  U \otimes_{A}R \cong U \propto KU \rightarrow 0.$$
It follows that $U$ is finitely presented and $H =U \propto KU$.

The inductive step is carried out just as we did for the case $n=1$
above, provided one substitutes the induction hypothesis for Lemma~\ref{3.2}. \\
$ii)\Longrightarrow i)$ Straightforward. $\Box$
\bigskip

{\noindent \it Proof of Theorem~\ref{3.1}.} (1) is already handled
by Theorem~\ref{2.8}(1). Specifically, $R(0,1)=0 \propto A$ is a
finitely generated ideal of $R$ which is not finitely presented
(by Lemma~\ref{3.3}).\\
(2) and (3) follow readily from a combination of Lemma~\ref{3.3}
with the next three facts:\\
(a) $R$ is a faithfully flat $A$-module.\\
(b) For $n\geq2$, a ring $B$ is $n$-coherent if and only if every
$(n-1)$-presented submodule of a finitely
generated free $B$-module is $n$-presented.\\
(c) For $n\geq2$ and $d\geq1$, a ring $B$ is an $(n,d)$-ring if
and only if every $(n-1)$-presented submodule of a finitely
generated free $B$-module has projective dimension $\leq d-1$.
$\Box$\bigskip

For $n\leq1$ or $d=0$, the $(n,d)$-property may not survive, in
general, in the trivial extension $R$ (even under strong
assumption on $A$). This is illustrated by the next example.

\begin{exmp}\label{3.4}
\rm Let $A$ be any arbitrary Pr\"ufer domain (i.e., $(1,1)$-domain) and let $R$
be the trivial ring extension of $A$ by its quotient field.
 Then $R$ is a $(2,1)$-ring which is neither a semihereditary ring
 (i.e., $(1,1)$-ring) nor a 2-von Neumann
regular ring (i.e., $(2,0)$-ring).
\end{exmp}
\proof   That $R$ is a $(2,1)$-ring which is not a $(1,1)$-ring is
ensured by Theorem~\ref{3.1}(3)\&(1), respectively. Now, let $J$
:=$R(a,0)$, where $a$ is a non-zero non-invertible element of $A$.
Since $(a,0)$ is a regular element of $R$, then the ideal $J$ of
$R$ has no non-zero annihilator. By \cite[Theorem 2.1]{M1}, $R$ is
not a 2-von Neumann regular ring. $\Box$
\bigskip

The B\'ezout property, however, does transfer reciprocally from
$A$ to $R$, as shown by the next result.

\begin{prop}\label{3.5}
  Let $R$ be as in Theorem~\ref{3.1}. Then $R$ is a B\'ezout ring if and
only if $A$  is a B\'ezout domain.
\end{prop}
\proof  Assume $R$ is a B\'ezout ring and let $I$ be a finitely
generated proper ideal of $A$. Then $J :=I \propto IK =I \propto
K$ is a finitely generated ideal of $R$. So $J =R(a,e)$ for some
$a \in A$ and $e \in E$. Therefore, $I =Aa$ and hence $A$ is a
B\'ezout domain.

Conversely, let $J$ be a finitely generated proper ideal of $R$.
Set $I :=\{a \in A / (a,e) \in J$ for some $e \in K\}$. We
consider two cases. {\sc Case 1}: $I =0$.  Necessarily, $J =0
\propto (1/b)L$ for some $b\not=0 \in A$ and some finitely
generated proper ideal $L$ of $A$. Further, $L=Aa$ since $A$ is a
B\'ezout domain. Hence $J =0 \propto A(a/b) =R(0,a/b)$, as
desired. {\sc Case 2}: $I \not= 0$. Let $(a,e) \in J$ such that $a
\not= 0$. Then, $(a,e)(0 \propto K) =0 \propto K \subseteq J$;
equivalently, $J  =I \propto IK =I \propto K$. But $I =Aa$ for
some $a \in A$ since $A$ is a B\'ezout domain. Therefore, $J =I
\propto K =R(a,0)$, completing the proof.
 $\Box$\bigskip

 Noteworthy is that new families of examples of non-semihereditary
 B\'ezout rings stem from the combination of Example~\ref{3.4} and
 Proposition~\ref{3.5}.\bigskip

 At this point, for the convenience of the reader, we recall from \cite{KM} the
 main result that establishes the transfer of the $(n,d)$-property
 to trivial ring extensions of local rings by their respective residue
 fields.

\begin{thm}\label{3.6} {\rm \cite[Theorem 1.1]{KM}}
Let $(A,M)$ be a local ring and let $R :=A \propto A/M$  be the
trivial ring extension of $A$  by
$A/M.$  Then: \\
{\rm(1)} $R$ is a $(3,0)$-ring provided $M$ is not finitely generated.\\
{\rm(2)} $R$ is not a $(2,d)$-ring, for each integer $d \geq 0$,
provided $M$ contains a regular element. $\Box$
\end{thm}

Clearly, Theorems \ref{3.1} and \ref{3.6} generate original
examples of $n$-coherent rings which, moreover, reflect no obvious
correlation between (strong) $n$-coherence and the large class of
finite conductor rings.

\begin{exmp}\label{3.7}\rm
Let $\Z$ be the ring of integers and $\Q =\qf(\Z)$. Then $R :=\Z
\propto \Q$ is a strong 2-coherent ring which is not a finite
conductor ring.
\end{exmp}
\proof Straightforward via Theorem~\ref{3.1} and
Theorem~\ref{2.8}. $\Box$

\begin{exmp} \label{3.8} \rm
 Let $(V,M)$ be a nondiscrete valuation domain.  Then $R
:=V \propto V/M$ is a 3-coherent ring which is neither 2-coherent
nor a finite conductor ring.
\end{exmp}
\proof $R$ is a $(3,0)$-ring by Theorem~\ref{3.6}(1) (since $M$ is
not a finitely generated ideal of $A$). Hence $R$ is $3$-coherent.
Further $R$ is not a $(2,0)$-ring by Theorem~\ref{3.6}(2). So by
\cite[Theorem 2.4]{C} $R$ is not $2$-coherent. Finally,
Theorem~\ref{2.6}(2) ensures that $R$ is not a finite conductor
ring. $\Box$
\end{subsection}
\begin{subsection}{Costa's Second Conjecture}

A well-known fact about semihereditary rings is that a ring $R$ is
a $(1,1)$-ring if and only if it is a weak $(1,1)$-ring. In this
vein, Costa's second conjecture is that the $(n,d)$-property and
the weak $(n,d)$-property are equivalent for any non-negative
integers $n$ and $d$. So far, it has been shown that this
conjecture is valid under the coherence assumption \cite{CK}. It
remains however elusively open, in general.\bigskip

Our modest objective in this subsection is mainly to test its
validity beyond the class of coherent rings. In this line, two
results are stated generating two new contexts of validity for
this conjecture. The first of these involves trivial ring
extensions issued from coherent domains.

\begin{thm}\label{3.9}
Let $A$ be a non-trivial coherent domain, $K =\qf(A)$, and $R:=A
\propto K$. Then, for any integers $n\geq2$ and $d\geq1$, $R$ is a
non-coherent ring such that the following statements are equivalent:  \\
{\rm (i)} $R$ is an $(n,d)$-ring;   \\
{\rm (ii)} $R$ is a weak $(n,d)$-ring;   \\
{\rm (iii)} $\w.dim(A) \leq d$. \end{thm}
\proof (i)$\Longleftrightarrow$(ii) is a straightforward
application of Theorem~\ref{3.1}(1)\&(3) combined to
\cite[Proposition 2]{CK}.\\
(i)$\Longrightarrow$(iii) Assume $R$ is an $(n,d)$-ring. By
Theorem~\ref{3.1}(3) and \cite[Theorem 2.4]{C}, $A$ is a
$(1,d)$-domain, hence
$\w.dim(A) \leq d$ by \cite[Theorem 1.3.9]{Gz2}. \\
(iii)$\Longrightarrow$(i) Assume that $\w.dim(A) \leq d$. Let $M
:=I \propto K$ be any arbitrary maximal ideal of $R$, where $I$ is
a maximal ideal of $A$ (Cf. \cite[Theorem 25.1]{H}). We have
$\fd_{A}(I) \leq d-1$ and so $\fd_{R}(M)\leq d-1$ (since $M \cong
I \otimes R$ and $R$ is $A$-flat). By \cite[Theorem 4.1]{C}, $R$
is a $(d+1,d)$-ring. Further, by Theorem~\ref{3.1}(2), $R$ is
strong $2$-coherent. It follows that $R$ is a $(2,d)$-ring (by
\cite[Theorem 2.4]{C}) and hence an $(n,d)$-ring, as desired.
$\Box$\bigskip

Recall that a ring $R$ is finite conductor if any ideal $I$ of $R$
with $\mu(I)\leq2$ is finitely presented, where $\mu(I)$ denotes
the cardinality of a minimal set of generators of $I$
\cite[Proposition 2.1]{Gz1}. Our next (and last) theorem tests
Costa's second conjecture in the class of finite conductor rings.
As might be expected, the ``$\mu(I)\leq2$" assumption (in the
above definition) restricts the scope of this result to $n=2$ and
$d\leq1$.

\begin{thm}\label{3.10}
Let $R$ be a finite conductor ring. Then:\\
{\rm (1)} The following statements are equivalent: \\
\hspace*{.2in} {\rm (i)} $R$ is a Pr\"ufer domain; \\
\hspace*{.2in} {\rm (ii)} $R$ is a $(2,1)$-domain; \\
\hspace*{.2in} {\rm (iii)} $R$ is a weak $(2,1)$-domain.\\
{\rm (2)} The following statements are equivalent: \\
\hspace*{.2in} {\rm (i)} $R$ is a von Neumann regular ring; \\
\hspace*{.2in} {\rm (ii)} $R$ is a $(2,0)$-ring;\\
\hspace*{.2in} {\rm (iii)} $R$ is a weak $(2,0)$-ring.
\end{thm}
\proof  (1) We need only prove $(iii)\Longrightarrow (i)$. Suppose
that (iii) holds. Let $I$ be an arbitrary $2$-generated ideal of
$R$ (i.e., $\mu(I)\leq2$). Then $I$ is finitely presented, and
hence projective by (iii). Therefore
$R$ is a Pr\"ufer domain by \cite[Theorem 22.1]{Gi}.\\
(2) We need only prove $(iii)\Longrightarrow (i)$. Suppose that
(iii) holds. It suffices to show that each principal ideal of $R$
is a direct summand of $R$. Let $I$ be a principal ideal of $R$.
Then $I$ is finitely presented by hypothesis, so that $R/I$ is a
$2$-presented cyclic $R$-module, and hence projective by (iii).
Therefore the following exact sequence splits:
$$0 \rightarrow  I \rightarrow  R \rightarrow  R/I \rightarrow 0$$
leading to the conclusion. $\Box$

\begin{rem} \rm
Assertion (1) of Theorem~\ref{3.10} cannot extend to rings with
zerodivisors. Indeed, let $R:=\Z \times E$ as in
Example~\ref{2.3}. Then $R$ is a finite conductor ring which is
not a semihereditary ring (since it is not coherent). On the other
hand, by \cite[Theorem 4.5]{C}, $R$ is a $(2,1)$-ring since
$wdim(R)=1$. $\Box$
\end{rem}
\end{subsection}
\end{section}

\bigskip\bigskip
\end{document}